\documentclass[12pt]{article}

\usepackage{theorem,amssymb,amsmath}
\advance \topmargin by -\headheight
\advance \topmargin by -\headsep
\evensidemargin \oddsidemargin
\author{J.-P. Allouche \\
CNRS, Institut de Math\'ematiques de Jussieu \\
\'Equipe Combinatoire et Optimisation \\
Universit\'e Pierre et Marie Curie, Case 247, 4 Place Jussieu \\
F-75252 Paris Cedex 05 \\
France \\
{\tt allouche@math.jussieu.fr}
\and
M. Mend\`es France \\
Universit\'e Bordeaux I \\
Math\'ematiques \\
F-33405 Talence Cedex \\
France \\
{\tt michel.mendes-france@math.u-bordeaux1.fr} \\
}

\title{Lacunary formal power series and the Stern-Brocot sequence}

\date{ }

\def \proof{\bigbreak\noindent{\it Proof.\ \ }}

\def \endpf{{\ \ $\Box$ \medbreak}}

\newtheorem{theorem}{Theorem}
\newtheorem{lemma}{Lemma}
\newtheorem{corollary}{Corollary}

\theorembodyfont{\rm}

\newtheorem{remark}{Remark}

\newtheorem{definition}{Definition}

\begin{document}

\maketitle

\begin{center}
{\it \`A la m\'emoire de Philippe Flajolet}
\end{center}

\begin{abstract}
Let $F(X) = \sum_{n \geq 0} (-1)^{\varepsilon_n} X^{-\lambda_n}$ be a real lacunary
formal power series, where $\varepsilon_n = 0, 1$  and $\lambda_{n+1}/\lambda_n > 2$.
It is known that the denominators $Q_n(X)$ of the convergents of its continued fraction
expansion are polynomials with coefficients $0, \pm 1$, and that the number of nonzero
terms in $Q_n(X)$ is the $n$th term of the Stern-Brocot sequence. We show that replacing
the index $n$ by any $2$-adic integer $\omega$ makes sense. We prove that $Q_{\omega}(X)$
is a polynomial if and only if $\omega \in {\mathbb Z}$. In all the other cases
$Q_{\omega}(X)$ is an infinite formal power series, the algebraic properties of which we
discuss in the special case $\lambda_n = 2^{n+1} - 1$.
\end{abstract}


\section{Introduction}

\subsection{Lacunary power series and continued fraction expansions}

Let $\Lambda = (\lambda_n)_{n \geq 0}$ be a sequence of positive integers with 
$0 < \lambda_0 < \lambda_1 < ...$ satisfying $\lambda_{n+1} / \lambda_n > 2$ for all 
$n \geq 0$. Consider the formal power series 
$F(X):= \sum_{n \geq 0} (-1)^{\varepsilon_n} X^{-\lambda_n}$,
where $\varepsilon_n = 0, 1$. As is well known, power series in $X^{-1}$ can be
represented by a continued fraction $[A_0(X), A_1(X), A_2(X), ...]$, where the
$A_j$'s are polynomials in $X$, and for all $i > 0$, $A_i(X)$ is a non-constant 
polynomial. Quite obviously, in the case of the above $F(X)$, one 
has $A_0(X) = 0$.

\bigskip

Let $P_n(X)/Q_n(X) = [0, A_1(X), A_2(X), ..., A_n(X)]$ be the $n$th convergent of 
$F(X)$. As was already discovered in \cite{ALMPS} and \cite{MPS}, the 
denominators $Q_n(X)$ are particularly interesting to study: their coefficients are 
$0, \pm 1$.

\subsection{A sequence of polynomials and a sequence of integers}

The denominators $Q_n(X)$ introduced above can be quite explicitly expressed (see \cite{MPS}):
$$
Q_n(X) = \sum_{k \geq 0} \sigma(k, \varepsilon) {\frac{n+k}{2} \choose k}_2 X^{\mu(k, \Lambda)}.
$$
The exponent of $X$ is given by $\mu(k, \Lambda) =  
\sum_{q \geq 0} e_q(k) (\lambda_q - \lambda_{q-1})$, with $\lambda_{-1} = 0$, where
$e_q(k)$ is the $q$th binary digit of $k = \sum_{q \geq 0} e_q(k) 2^q$. The sign
of the monomials is given by $\sigma(k, \varepsilon) = (-1)^{\nu(k) + \bar{\mu}(k, \varepsilon)}$
where $\nu(k)$ is the number of occurrences of the block $10$ in the usual left-to-right reading
of the binary expansion of $k$ (e.g., $\nu(\mbox{\rm twelve}) = 1$), and where 
$\bar{\mu}(k, \varepsilon) = \sum_{q \geq 0} e_q(k) (\varepsilon_{k-1} - \varepsilon_{k-2})$,
with $\varepsilon_{-1} = \varepsilon_{-2} = 0$. The symbol ${a \choose b}_2$ is the {\bf integer} 
equal to $0$ or $1$, according to the value modulo $2$ of the binomial coefficient ${a \choose b}$, 
with the following convention: if $a$ is not an integer, or if $a$ is a positive integer and 
$a < b$, then ${a \choose b} := 0$. For example, as soon as $n$ and $k$ have opposite parities, 
${\frac{n+k}{2} \choose k}_2 = 0$. In \cite{ALMPS} it was observed that
the number of non-zero monomials in $Q_n(X)$ is $u_n$, the $n$th term of the celebrated
Stern-Brocot sequence defined by $u_0 = u_1 = 1$, and the recursive relations
$u_{2n} = u_n + u_{n-1}$, $u_{2n+1} = u_n$ for all $n \geq 1$.
This sequence is also called the Stern diatomic series (see sequence A002487 in 
\cite{OEIS}). It was studied by several authors, see, e.g., \cite{DS1} and its list of 
references (including the historical references \cite{Brocot, Stern}), see also 
\cite{Urbiha, Northshield}, or see \cite{HKMPP} for a relation 
between the Stern sequence and the Towers of Hanoi. (Note that some authors have the 
slightly different definition: $v_0 = 0$, $v_{2n} = v_n$, $v_{2n+1} = v_n + v_{n+1}$; 
clearly $u_n = v_{n+1}$ for all $n \geq 0$.)

\bigskip

Our purpose here is to pursue our previous discussions on the sequence of polynomials
$Q_n(X)$ in relationship with the Stern-Brocot sequence. 

\begin{remark}
The sequence $(\nu(n))_{n \geq 0}$ happens to be related to the paperfolding sequence.
Indeed, define $v(n) := (-1)^{\nu(n)}$ and $w(n) := v(n)v(n+1)$. From the definition 
of $\nu$, we have for every $n \geq 0$ the relations $v(2n+1) = v(n)$, $v(4n) = v(2n)$,
and $v(4n+2) = - v(n)$. Equivalently, for every $n \geq 0$, we have $v(2n+1) = v(n)$, and
$v(2n) = (-1)^n v(n)$. Hence, for every $n \geq 0$, we have $w(n) = v(2n)v(2n+1) = 
(-1)^n (v(n)^2) = (-1)^n$, and $w(2n+1) = v(2n+1)v(2n+2) = (-1)^{n+1} v(n) v(n+1) = 
(-1)^{n+1} w(n)$. It it then clear that, if $z(n) := w(2n+1)$, then $z(2n) = -w(2n) = -(-1)^n$
and $z(2n+1) = z(n)$. In other words the sequence $(z(n))_{n \geq 0}$ is the classical
paperfolding sequence, and the sequence $(w(n))_{n \geq 0}$ itself is a paperfolding sequence,
see e.g., \cite[p.~125]{mmfvdp} where the sequences are indexed by $n \geq 1$ instead of
$n \geq 0$.
\end{remark}

\subsection{A partial order on the integers}\label{order}

Let $m = e_0(m) e_1(m)...$ and $k = e_0(k) e_1(k)...$ be two nonnegative integers
together with their binary expansion, which of course terminates with a tail of $0$'s.
Lucas \cite{Lucas} observed that
$$
{m \choose k} \equiv \prod_{i \geq 0} {e_i(m) \choose e_i(k)} \bmod 2.
$$
This implies the following relation (in ${\mathbb Z}$)
$$
{m \choose k}_2 = \prod_{i \geq 0} {e_i(m) \choose e_i(k)},
$$
so that we have ${m \choose k}_2 = 1$ if and only if $e_i(k) \leq e_i(m)$ 
for all $i \geq 0$.

\noindent
We will say that $m$ {\it dominates} $k$ and we write $k <\!\!< m$, if 
$e_i(k) \leq e_i(m)$ for all $i \geq 0$. In other words the sequence 
$k \to {m \choose k}_2$ is the characteristic function of the $k$'s dominated by $m$.
(This order was used in, e.g., \cite{AFM}.)

\bigskip

As a consequence of our remarks, the Stern-Brocot sequence has the following 
representation
$$
u_n = \sum_{k < \!\!< \frac{k+n}{2}} 1.
$$

\begin{remark}\label{carlitz}
This last relation can be easily deduced from a result of Carlitz 
\cite{Carlitz62, Carlitz64} (Carlitz calls $\theta_0(n)$ what we call $u_n$):
$$
u_n = \sum_{0 \leq 2r \leq n} {n-r \choose r}_2.
$$
Indeed, we have
$$
\begin{array}{lll}
\displaystyle\sum_{k < \!\!< \frac{k+n}{2}} 1 &=& 
\displaystyle\sum_{\stackrel{\scriptstyle 0 \leq k \leq n}{k \equiv n \bmod 2}}
{\frac{k+n}{2} \choose k}_2
= \sum_{\stackrel{\scriptstyle 0 \leq k' \leq n}{k' \equiv 0 \bmod 2}} 
{n-\frac{k'}{2} \choose n-k'}_2 \ \ \ \mbox{\rm (by letting $k' = n-k$)} \\
&=& \displaystyle\sum_{0 \leq 2r \leq n} {n-r \choose n-2r}_2
= \sum_{0 \leq 2r \leq n} {n-r \choose r}_2 \ \ \
\mbox{\rm (by using ${a \choose b} = {a \choose a-b}$)}.
\end{array}
$$
Also note that in \cite{Carlitz64} the range $0 \leq 2r < n$ should be replaced
by $0 \leq 2r \leq n$ as in \cite{Carlitz62} (see also \cite[Corollary~6.2]{DS1} 
where the index $n$ should be adjusted). Let us finally indicate that this remark is 
also Corollary~13 in \cite{ALMPS}.
\end{remark}

\begin{remark}
The relation $u_n = \sum_{0 \leq 2r \leq n} {n-r \choose r}_2$ can give the idea 
(inspired by the classical {\em binomial transform}) of introducing a map on sequences 
$(a_n)_{n \geq 0} \to (b_n)_{n \geq 0}$ with 
$b_n := \sum_{0 \leq 2r \leq n} {n-r \choose r}_2 a_r$, so that in particular the
image of the constant sequence $1$ is the Stern-Brocot sequence. One can also go a 
step further by defining a map ${\mathcal C}$ which associates with two sequences 
${\mathbf a} = (a_n)_{n \geq 0}$ and ${\mathbf b} = (b_n)_{n \geq 0}$ the sequence
$$
{\mathcal C}({\mathbf a}, {\mathbf b}) := 
\left(\sum_{0 \leq 2r \leq n} {n-r \choose r}_2 a_r b_{n-r}\right)_{n \geq 0}.
$$
It is unexpected that some variations on the Stern-Brocot sequences (different from but
in the spirit of the twisted Stern sequence of \cite{Bacher}) are related to the celebrated 
Thue-Morse sequence (see, e.g., \cite{AS-ubiq}). In fact, recall that the $\pm 1$ Thue-Morse 
sequence ${\mathbf t} = ((t_n)_{n \geq 0}$ can be defined by $t_0 = 1$ and, for all $n \geq 0$,
$t_{2n}= t_n$ and $t_{2n+1} = -t_n$. Now define the sequences
${\mathbf \alpha} = (\alpha_n)_{n \geq 0}$, ${\mathbf \beta} = (\beta_n)_{n \geq 0}$, 
${\mathbf \gamma} = (\gamma_n)_{n \geq 0}$ by
$$
{\mathbf \alpha} := {\mathcal C}({\mathbf t}, {\mathbf 1}), \
{\mathbf \beta}  := {\mathcal C}({\mathbf 1}, {\mathbf t}), \
{\mathbf \gamma} := {\mathcal C}({\mathbf t}, {\mathbf t}). 
$$
Then the reader can check that these sequences satisfy respectively
$$
\begin{array}{llll}
&{\mathbf \alpha}(0) = 1, \ {\mathbf \alpha}(1) = 1, \ &\mbox{\rm and for all $n \geq 1$}, \
&\alpha_{2n} = \alpha_n - \alpha_{n-1}, \ \alpha_{2n+1} = \alpha_n \\
&{\mathbf \beta}(0) = 1, \ {\mathbf \beta}(1) = -1, \ &\mbox{\rm and for all $n \geq 1$}, \
&\beta_{2n} = \beta_n - \beta_{n-1}, \ \beta_{2n+1} = -\beta_n \\
&{\mathbf \gamma}(0) = 1, \ {\mathbf \gamma}(1) = -1, \ &\mbox{\rm and for all $n \geq 1$}, \
&\gamma_{2n} = \gamma_n + \gamma_{n-1}, \ \gamma_{2n+1} = -\gamma_n \\
\end{array}
$$
so that, with the notation of \cite{OEIS},
$$
(\alpha_n)_{n \geq 0} = (A005590(n+1))_{n \geq 0}, \ 
(\beta_n)_{n \geq 0} = (A177219(n+1))_{n \geq 0}, \
(\gamma_n)_{n \geq 0} = (A049347(n))_{n \geq 0}.
$$
The last sequence $(\gamma_n)_{n \geq 0}$ is the $3$-periodic sequence with period 
$(1, -1, 0)$ (hint: prove by induction on $n$ that for all $j \leq n$ one has
$(\gamma_{3j}, \gamma_{3j+1},\gamma_{3j+2}) = (1, -1, 0)$).
\end{remark}

\section{More on the sequence $Q_n(X)$ and a note on $P_n(X)$ for a special $\Lambda$}

We now specialize to the case $\lambda_n = 2^{n+1} - 1$. In that case, $\mu(k, \Lambda) = k$.
Also note that $\sigma(k, \varepsilon) \equiv 1 \bmod 2$. Let $P_n(X)/Q_n(X)$ denote as 
previously the $n$th convergent of the continued fraction of the formal power series 
$\sum_{i \geq 1} (-1)^{\varepsilon_i} X^{1-2^i}$. We begin with a short section on $P_n$. 
The rest of the section will be devoted to the ``simpler'' polynomials $Q_n$.

\subsection{The sequence $P_n$ modulo $2$}

\begin{theorem} 
We have $P_n(X) \equiv Q_{n-1}(X) \bmod 2$ for $n \geq 1$.
\end{theorem}

\proof Let $F(X) = \sum_{i \geq 1} (-1)^{\varepsilon_i} X^{1-2^i}$. Define the formal 
power series $\Phi(X)$ by its continued fraction expansion $\Phi(X) = [0, X, X, \ldots]$. 
Its $n$th convergent is given by $\pi_n(X)/\kappa_n(X) = [0, X, \ldots, X]$ ($n$ partial 
quotients equal to $X$). An immediate induction shows that $\pi_n(X) = \kappa_{n-1}$ for 
$n \geq 1$.
Reducing $F(X)$ modulo $2$, we see that $F^2(X) + XF(X) + 1 \equiv 0 \bmod 2$. On the 
other hand $\Phi(X) = 1/(X + \Phi(X))$, hence $\Phi^2(X) + X\Phi(X) + 1 \equiv 0 \bmod 2$.
This implies that $F(X) \equiv \Phi(X) \bmod 2$. Hence $P_n(X) \equiv \pi_n(X) \bmod 2$
and $Q_n(X) \equiv \kappa_n(X) \bmod 2$: to be sure that the convergents of the reduction
modulo $2$ of $F$ are equal to the reduction modulo $2$ of the convergents of $F(X)$, the
reader can look at, e.g., \cite{vdP-reduc}. Thus $P_n(X) \equiv \pi_n(X) = \kappa_{n-1}(X)
\equiv Q_{n-1}(X) \bmod 2$.

\begin{corollary}\label{rel}
The following congruence is satisfied by $Q_n(X)$ for $n \geq 1$:
$$
Q_n^2(X) - Q_{n+1}(X)Q_{n-1}(X) \equiv 1 \bmod 2.
$$
\end{corollary}

\proof Use the classical identity $P_{n+1}(X)Q_n(X) - P_n(X)Q_{n+1}(X)  = (-1)^n$
for the convergents of a continued fraction.

\subsection{The sequence $Q_n$ and the Chebyshev polynomials}

We have the formula
$$
Q_n(X) \equiv \sum_{k \geq 0} {\frac{n+k}{2} \choose k}_2 X^k
\equiv \sum_{\stackrel{\scriptstyle 0 \leq k \leq n}{k \equiv n \bmod 2}}
{\frac{k+n}{2} \choose k}_2 X^k \ \ \ \ \bmod 2.
$$
The Chebyshev polynomials of the second kind (see, e.g., 
\cite[p.~184--185]{Erdelyi}) are defined by
$$ 
U_n(\cos \theta) = \frac{\sin (n+1)\theta}{\sin \theta}\cdot
$$
They have the well-known explicit expansion
$$
U_n(X) = \sum_{0 \leq k \leq n/2} (-1)^k {n-k \choose k} (2X)^{n-2k}.
$$
We thus get a relationship between $Q_n$ and $U_n$ (compare with the related but not
identical result \cite[Proposition~6.1]{DS1}).

\begin{theorem}\label{Cheb}
The reductions modulo $2$ of $Q_n(X)$ and of $U_n(X/2)$ are equal.
\end{theorem}

\proof We can write modulo $2$
$$
\begin{array}{lll}
Q_n(X) 
&\equiv& \displaystyle\sum_{\stackrel{\scriptstyle 0 \leq k' \leq n}{k' \equiv 0 \bmod 2}}
{n-\frac{k'}{2} \choose n-k'}_2 X^{n-k'} \ \ \ 
\mbox{\rm (by letting $k' = n-k$)} \\
&\equiv& \displaystyle\sum_{0 \leq 2r \leq n} {n-r \choose n-2r}_2 X^{n-2r}
\\
&\equiv& \displaystyle\sum_{0 \leq 2r \leq n} {n-r \choose r}_2 X^{n-2r}
\ \ \ \mbox{\rm (by using ${a \choose b} = {a \choose a-b}$)}.
\end{array}
$$
Hence $Q_n(X) \equiv U_n(X/2) \bmod 2$. \endpf

\bigskip

As an immediate application of Theorem~\ref{Cheb} (and of Remark~\ref{carlitz}) 
we have the following corollaries.

\begin{corollary}
The number of odd coefficients in the (scaled) Chebyshev polynomial of the second kind
$U_n(X/2)$ is equal to the Stern-Brocot sequence $u_n$.
\end{corollary}

\begin{remark}
Corollary~\ref{rel} above can also be deduced from Theorem~\ref{Cheb} using a classical 
relation for Chebyshev polynomials implied by their expression using sines.
\end{remark}

\begin{remark}
The polynomials $Q_n(X)$ are also related to the Fibonacci polynomials (see, e.g.,
\cite{Ency}) and to Morgan-Voyce polynomials which are a variation on the Chebyshev 
polynomials (on Morgan-Voyce polynomials, introduced by Morgan-Voyce in dealing with 
electrical networks see, e.g., \cite{Swamy, AJ, GR} and the references therein).
Indeed, the Fibonacci polynomials satisfy
$$
F_{n+1}(X) = \sum_{2j \leq n} {n-j \choose j} X^{n-2j}
$$
(compare with the proof of Theorem~\ref{Cheb}),
while the Morgan-Voyce polynomials satisfy
$$
b_n(X) = \sum_{k \leq n} {n+k \choose n-k} X^k \ \ \mbox{\rm and} \ \
B_n(X) = \sum_{k \leq n} {n+k+1 \choose n-k} X^k
$$
(note that ${n+k \choose n-k} = {n+k \choose 2k}$, that
${n+k+1 \choose n-k} = {n+k+1 \choose 2k+1}$, and see 
Lemmas~\ref{digit} and \ref{kernel} below).
\end{remark}

\begin{remark}
The polynomials that we have defined are related to the Stern-Brocot sequence, but they 
differ from Stern polynomials occurring in the literature, in particular they are not the 
same as those introduced in \cite{KMP}. They also differ from the polynomials studied in 
\cite{DS1, DS2}.
\end{remark}

\subsection{Extension of $Q_n(X)$ to $Q_{\omega}(X)$ with $\omega \in {\mathbb Z}_2$}
\label{extension}

\begin{definition}\label{binom}
Let $\omega = \sum_{i \geq 0} \omega_i 2^i = \omega_0 \omega_1 \omega_2 
\ldots \in {\mathbb Z}_2$ be a $2$-adic integer, or equivalently an 
infinite sequence of $0$'s and $1$'s. For a nonnegative integer $k$
whose binary expansion is given by $k = \sum_{i \geq 0} k_i 2^i$, we define
$$
{\omega \choose k}_2 = \prod_{i \geq 0} {\omega_i \choose k_i}.
$$
\end{definition}

The infinite product ${\omega \choose k}_2$ is well defined since, for large $i$, 
${\omega_i \choose k_i}$ reduces to ${\omega_i \choose 0} =  1$. It is equal to $0$ or 
$1$. The above product extends Lucas' observation to all $2$-adic integers $\omega$. 
In particular, since $-1 = \sum_{i \geq 0} 2^i = 1^{\infty}$, we see that $-1$ dominates 
all $k \in {\mathbb N}$ (where the order introduced in Section~\ref{order} is generalized 
in the obvious way). A similar definition (binomials and order) occurs in \cite{MFvdP}.

\begin{definition}
In the general case for $\Lambda$, with $\lambda_{n+1}/\lambda_n > 2$, and $\varepsilon = 0, 1$,
the polynomials $Q_n(X)$ above naturally extend to formal power series $Q_{\omega}(X)$ 
defined for $\omega = \omega_0 \omega_1 \omega_2 \ldots \in {\mathbb Z}_2$ by
$$
Q_{\omega}(X) = 
\sum_{k \geq 0} \sigma(k, \varepsilon) {\frac{\omega+k}{2} \choose k}_2 X^{\mu(k, \Lambda)} 
= \sum_{\stackrel{\scriptstyle k \equiv \omega \bmod 2}{k < \!\! < \frac{\omega+k}{2}}} 
\sigma(k, \varepsilon) X^{\mu(k, \Lambda)}.
$$
\end{definition}

\begin{remark}\label{negative}
The reader can check (e.g., by using integer truncations of $\omega$ tending
to $\omega$) that
$$
{\omega \choose k} \equiv {\omega \choose k}_2 \bmod 2
$$
where the binomial coefficient ${\omega \choose k}$ is defined by
$$
{\omega \choose k} = \frac{\omega (\omega-1)\ldots (\omega-k+1)}{k!} \in {\mathbb Z}_2.
$$
In particular, we see that for any $2$-adic integer $\ell$,
$$
{- \ell \choose k} = (-1)^k {\ell + k -1 \choose k}, \ \ \mbox{\rm hence} \ \
{- \ell \choose k}_2 = {\ell + k -1 \choose k}_2.
$$
Now for $n \in {\mathbb N}$ we have 
$$
Q_{-n}(X) = \sum_{k \geq 0} \sigma(k, \varepsilon) {\frac{-n+k}{2} \choose k}_2 X^{\mu(k, \Lambda)}
= \sum_{k \geq 0} \sigma(k, \varepsilon) {-\frac{(n-k)}{2} \choose k}_2 X^{\mu(k, \Lambda)},
$$
thus 
$$
Q_{-n}(X) = \sum_{k \geq 0} \sigma(k, \varepsilon) {\frac{n-k}{2}+k-1 \choose k}_2 X^{\mu(k, \Lambda)}
= \sum_{k \geq 0} \sigma(k, \varepsilon) {\frac{n-2+k}{2} \choose k}_2 X^{\mu(k, \Lambda)} 
= Q_{n-2}(X).
$$
In particular $Q_{-n}$ and $Q_{n-2}$ have same degree. Also note that the definition
of $Q_{-n}$ for $n \in {\mathbb N}$ yields
$$
Q_{-1}(X) = \sum_{k \geq 0} \sigma(k, \varepsilon) {\frac{k-1}{2} \choose k}_2 X^{\mu(k, \Lambda)} 
= 0.
$$
\end{remark}

\begin{remark}
If $\lambda_n = 2^{n+1} - 1$, Corollary~\ref{rel} can be extended to $2$-adic integers: using 
again truncations of $\omega$ tending to $\omega$ yields, for any $2$-adic integer $\omega$,
$$
Q_{\omega}^2(X) - Q_{\omega+1}(X)Q_{\omega-1}(X) \equiv 1 \bmod 2.
$$
\end{remark}

\subsection{Extension of the sequence $(u_n)_{n \geq 0}$ to negative indices}

What precedes suggests two ways of extending the sequence $(u_n)_{n \geq 0}$ to 
negative integer indices. First, we noted the relation 
$u_n = \displaystyle\sum_{k < \! \! < \frac{n+k}{2}} 1$, i.e.,
$u_n$ is the number of monomials with non-zero coefficient in $Q_n(X)$.
But from the previous section, we can define $Q_{-n}(X)$ for $n \in {\mathbb N}$,
and we have that $Q_{-n}(X) = Q_{n-2}(X)$. This suggests the definition
$$
u_{-n} := u_{n-2} \ \mbox{\rm for all } n \geq 2.
$$
Strictly speaking, this definition leaves the value $u_{-1}$ indeterminate,
but, since $u_n$ is the number of monomials with nonzero coefficients in $Q_n$,
the remark above that $Q_{-1} = 0$ implies $u_{-1} = 0$.

\bigskip

Another way of generalizing $u_n$ to negative indices would be to
use the recursion 
$$
u_{2n} = u_n + u_{n-1}, \ \ u_{2n+1} = u_n, \ \ \mbox{\rm for all } n \geq 1,
$$
allowing non-positive values for $n$. Allowing first $n=0$ leads to
$u_0 = u_0 + u_{-1}$, hence $u_{-1} = 0$. On the other hand we claim that
the relation $u_{-n} := u_{n-2}$ for all $n \geq 2$ leads to the same
recursion formulas with $u_{2n}$ and $u_{2n+1}$ for non-positive $n$. 
Indeed, let $m = -n$, with $n \geq 2$. Then
$$
u_{2m} = u_{-2n} = u_{2n-2} = u_{2(n-1)} = u_{n-1} + u_{n-2} 
= u_{-n-1} + u_{-n} = u_{m-1} + u_m
$$
and
$$
u_{2m+1} = u_{-2n+1} = u_{2n-3} = u_{2(n-2)+1} = u_{n-2} = u_{-n} = u_m.
$$
We thus finally have a generalization compatible with both approaches, yielding
$$
\ldots \ u_{-4} = 2, \ u_{-3} = 1, \ u_{-2} = 1, \ u_{-1} = 0, \ u_0 = 1, 
\ u_1 = 1, \ u_2 = 2, \ u_3 = 1, \ u_4 = 3 \ \ldots
$$
and the definition
\begin{definition}
The Stern-Brocot sequence $(u_n)_{n \geq 0}$ can be extended to a sequence
$(u_n)_{n \in {\mathbb Z}}$ by letting $u_{-n} = u_{n-2}$ for $n \geq 2$, and
$u_{-1} = 0$. This sequence satisfies the same recursive relations as the
initial sequence $(u_n)_{n \geq 0}$, namely $u_{2n} = u_n + u_{n-1}$ and
$u_{2n+1} = u_n$ for all $n \in {\mathbb Z}$.
\end{definition}

\section{The arithmetical nature of the power series $Q_{\omega}(X)$}

Recall that the formal series $Q_{\omega}(X)$, where $\omega = \omega_0 \omega_1 \ldots$
belongs to ${\mathbb Z}_2$, is given by
$$
Q_{\omega}(X) 
= \sum_{k \geq 0} \sigma(k, \varepsilon) {\frac{\omega+k}{2} \choose k}_2 X^{\mu(k, \Lambda)} 
= \sum_{\stackrel{\scriptstyle k \equiv \omega \bmod 2}{k < \!\! < \frac{\omega+k}{2}}} 
\sigma(k, \varepsilon) X^{\mu(k, \Lambda)}.
$$
We have seen that $Q_{\omega}(X)$ reduces to a polynomial if $\omega$ belongs to
${\mathbb Z}$. We will prove that this is a necessary and sufficient condition for
this series being a polynomial. Then we will address the question of the algebraicity
of $Q_{\omega}(X)$, on ${\mathbb Q}(X)$ and on ${\mathbb Z}/2{\mathbb Z}(X)$, in the 
special case $\lambda_n = 2^{n+1} - 1$. We begin with a lemma.

\begin{lemma}\label{digit}
Let $\omega = \omega_0 \omega_1 \ldots$ belong to ${\mathbb Z}_2$.
Then the following properties hold.

\begin{itemize}

\item[{\rm (i)}] For every $j \geq 0$,
$$
{\omega + 2^j \choose 2^{j+1}}_2 \equiv \omega_j + \omega_{j+1} \bmod 2.
$$

\item[{\rm (ii)}] The sequence $({\omega + 2^j \choose 2^{j+1}}_2)_{j \geq 0}$ is
ultimately periodic if and only if $\omega$ is rational.

\item[{\rm (iii)}] The sequence $({\omega + 2^j \choose 2^{j+1}}_2)_{j \geq 0}$ is
ultimately equal to $0$ if and only if $\omega$ is an integer.

\item[{\rm (iv)}] For every $k \geq 0$
$$
{\frac{\omega + k}{2} \choose k}_2 = {\omega + k + 1 \choose 2k+1}_2.
$$

\item[{\rm (v)}] If $\omega \neq -1$, there exist an integer $\ell \geq 0$ and 
a $2$-adic integer $\omega'$ such that $\omega = 2^{\ell}-1 + 2^{\ell +1} \omega'$.  
Let $f_\omega(k) := 
{\frac{\omega + k}{2} \choose k}_2 = {\omega + k+1 \choose 2k+1}_2$. Then for
any integer $k'$ we have $f_{\omega}(2^{\ell} - 1 + 2^{\ell + 1}k') =
{\omega' + k' \choose 2k'}_2$.

\item[{\rm (vi)}] Suppose that there exist $\ell \geq 0$ and $j \geq 0$ such that
$\omega = 2^{\ell} - 1 + 2^{\ell +1}(2^j(2\omega'+1))$. Then for any integer
$k'$ we have 
$f_{\omega}(2^{\ell} - 1 + 2^{\ell + 1}(2^j(2k'+1))) =
{\omega' + k' + 1 \choose 2k'+1}_2$.

\end{itemize}
\end{lemma}

\proof In order to prove (i) we write
$$
\begin{array}{llllllllll}
\omega + 2^j &=& & &\omega_0 &\omega_1 &\ldots &\omega_j &\omega_{j+1} &\ldots \\
&\ \ & &+ &0 &0 &\ldots &1 &0 &\ldots \\
&=& & &\omega_0 &\omega_1 &\ldots &\alpha_j &\alpha_{j+1} &\ldots \\
\end{array}
$$
where $\alpha_j$ and $\alpha_{j+1}$ are given by
$$
\begin{array}{lll}
&\mbox{\rm if $\omega_j = 0$ and $\omega_{j+1} = 0$,}  \
&\mbox{\rm then $\alpha_j = 1$ and $\alpha_{j+1} = 0$} \\
&\mbox{\rm if $\omega_j = 0$ and $\omega_{j+1} = 1$,}  \
&\mbox{\rm then $\alpha_j = 1$ and $\alpha_{j+1} = 1$} \\
&\mbox{\rm if $\omega_j = 1$ and $\omega_{j+1} = 0$,}  \
&\mbox{\rm then $\alpha_j = 0$ and $\alpha_{j+1} = 1$} \\
&\mbox{\rm if $\omega_j = 1$ and $\omega_{j+1} = 1$,}  \
&\mbox{\rm then $\alpha_j = 0$ and $\alpha_{j+1} = 0$}. \\
\end{array}
$$
By inspection we see that $\alpha_{j+1} \equiv \omega_j + \omega_{j+1} \bmod 2$.
Now we write
$$
{\omega + 2^j \choose 2^{j+1}}_2 
= \left(\prod_{0 \leq k \leq j-1} {\omega_k \choose 0}_2\right)
  {\alpha_j \choose 0}_2 {\alpha_{j+1} \choose 1}_2 
  \left(\prod_{k \geq j+2} {\alpha_k \choose 0}_2\right)
= \alpha_{j+1} \equiv \omega_j + \omega_{j+1} \bmod 2. 
$$

\bigskip

Let us prove (ii). We note that the sequence 
$((\omega_j + \omega_{j+1}) \bmod 2)_{j \geq 0}$ is ultimately periodic 
if and only if the sequence $(\omega_j \bmod 2)_{j \geq 0}$ is ultimately 
periodic (hence if and only if the sequence $(\omega_j)_{j \geq 0}$ itself is
ultimately periodic): indeed, $((\omega_j + \omega_{j+1}) \bmod 2)_{j \geq 0}$ is 
ultimately periodic if and only if the formal power series
$G(X) := \sum_{j \geq 0}(\omega_j + \omega_{j+1}) X^j$ is rational (as an element of
${\mathbb Z}/2{\mathbb Z}[[X]]$). But, if we let $H(X)$ denote the formal power
series $H(X) := \sum_{j \geq 0} \omega_j X^j \in {\mathbb Z}/2{\mathbb Z}[[X]]$, then
$XG(X)+\omega_0 = (1+X) H(X)$. So $G(X)$ is rational if and only if $H$ is, if and only
if $(\omega_j \bmod 2)_{j \geq 0}$ is ultimately periodic, i.e., if the $2$-adic integer
$\omega$ is rational.

\bigskip

To prove (iii), we note that ${\omega+2^j \choose 2^{j+1}}_2 = 0$ for $j$ large enough
implies from Lemma~\ref{digit} (i) that $\omega_j + \omega_{j+1} \equiv 0 \bmod 2$ for
$j$ large enough. This means that $\omega_j \equiv \omega_{j+1} \bmod 2$ for $j$
large enough, or equivalently $\omega_j = \omega_{j+1}$ for $j$ large enough.
But then either $\omega_j = \omega_{j+1} = 0$ for large $j$, hence $\omega$ is
a nonnegative integer, or $\omega_j = \omega_{j+1} = 1$ for large $j$, hence
$\omega$ is a negative integer. We thus finally get that $\omega$ belongs to
${\mathbb Z}$. The converse is straightforward.

\bigskip

We prove (iv) by considering the parities of $\omega$ and $k$.
First note that if $\omega$ and $k$ have opposite parities, then
${\frac{\omega + k}{2} \choose k}_2 = 0$ while ${\omega + k +1 \choose 2k+1}_2 = 0$
(use Definition~\ref{binom} and look at the last digit of $\omega+k+1$ and of $2k+1$).
Now if $\omega = 2\omega'$ and $k = 2k'$, we have
${\frac{\omega + k}{2} \choose k}_2 = {\omega' + k' \choose 2k'}_2$ while
${\omega + k +1 \choose 2k+1}_2 = {2(\omega' + k') + 1 \choose 4k' +1}_2
= {\omega' + k' \choose 2k'}_2$ (use Definition~\ref{binom} again).
Finally if $\omega = 2\omega'+1$ and $k = 2k'+1$, we have
${\frac{\omega + k}{2} \choose k}_2 = {\omega' + k' + 1 \choose 2k'+1}_2$ while
${\omega + k +1 \choose 2k+1}_2 = {2(\omega' + k' + 1) + 1 \choose 4k'+3}_2 =
{\omega' + k' + 1 \choose 2k'+1}_2$ (use Definition~\ref{binom} again).

\bigskip

Let us prove (v). Since $\omega \neq -1$, its $2$-adic expansion contains at least 
one zero. Write $\omega = 1 1 \ldots 1 0 \omega_{\ell+1} \omega_{\ell+2} \ldots$, so 
that the $2$-adic expansion of $\omega$ begins with exactly $\ell \geq 0$ ones. 
Defining $\omega' := \omega_{\ell+1} \omega_{\ell+2} \ldots$, we thus have
$\omega = 2^{\ell} - 1 + 2^{\ell + 1}\omega'$. Now for any integer $k'$ we have from
Definition~\ref{binom}
$$
f_{\omega}(2^{\ell} - 1 + 2^{\ell + 1}k') =
{\omega + 2^{\ell} + 2^{\ell + 1}k' \choose 2^{\ell + 1} - 1 + 2^{\ell + 1}(2k')}_2
= {2^{\ell+1} - 1 + 2^{\ell + 1}(\omega'+k') \choose 2^{\ell + 1} - 1 + 2^{\ell + 1}(2k')}_2
= {\omega' + k' \choose 2k'}_2.
$$

\bigskip

We finally prove (vi). Using (v) we see that
$$
f_{\omega}(2^{\ell} - 1 + 2^{\ell + 1}(2^j(2k'+1)))
= {2^j(2\omega'+1+2k'+1) + 1 \choose 2^{j+1}(2k'+1) + 1}_2
= {\omega'+k'+1 \choose 2k'+1}_2.
$$
\bigskip

Now we can prove the following result.

\begin{theorem}\label{polynomial}
Let $\omega$ be a $2$-adic integer. The formal power series $Q_{\omega}(X)$
is a polynomial if and only if $\omega$ belongs to ${\mathbb Z}$.
\end{theorem}

\proof If $n$ is a nonnegative integer, then $Q_n(X)$ is a polynomial. So is 
$Q_{-n}(X)$ for $n \neq 1$ because $Q_{-n} = Q_{n-2}$ as we have seen in Remark~\ref{negative}.
On the other hand $Q_{-1}(X)$ is also a polynomial since $Q_{-1}(X) = 0$. 
Conversely suppose that $Q_{\omega}(X)$ is a polynomial, for some 
$\omega = \omega_0 \omega_1 \ldots$ in ${\mathbb Z}_2$. The coefficients of the monomials
$X^{\mu(k, \Lambda)}$ in $Q_{\omega}(X)$, i.e., 
$\sigma(k, \varepsilon) {\frac{\omega+k}{2} \choose k}_2$, are equal to zero for $k$ large enough. 
Thus $f_{\omega}(k) = {\frac{\omega+k}{2} \choose k}_2$ is equal to zero for $k$ large enough. 
We may suppose that $\omega \neq -1$; thus, using the notation in Lemma~\ref{digit} (v), 
we certainly have that $f_{\omega}(2^{\ell} -1 + 2^{\ell+1} k') =0$ for $k'$ large enough.
Using Lemma~\ref{digit} (v), we thus have ${\omega'+k' \choose 2k'}_2 = 0$ for $k'$ large 
enough. This implies ${\omega'+2^j \choose 2^{j+1}}_2 = 0$ for $j$ large enough.
Lemma~\ref{digit} (iii) yields that $\omega'$, hence $\omega$, belongs to ${\mathbb Z}$. 

\bigskip

Before proving our Theorem~\ref{algebr} characterizing the algebraicity 
of the series $Q_{\omega}(X)$ for a special $\Lambda$, we need a lemma.

\begin{lemma}\label{kernel}

Let $\omega = \omega_0 \omega_1 \ldots$ be a $2$-adic integer. We let 
$(f_{\omega}(k))_{k \geq 0}$, $(g_{\omega}(k))_{k \geq 0}$, $(h_{\omega}(k))_{k \geq 0}$ 
denote the sequences
$$
f_{\omega}(k) := {\omega + k + 1 \choose 2k+1}_2, \ \
g_{\omega}(k) := {\omega + k \choose 2k}_2, \ \
h_{\omega}(k) := {\omega + k \choose 2k+1}_2.
$$
Then we have the following relations.
$$
\begin{array}{lllll}
& f_{2\omega}(2k) = g_{\omega}(k), & f_{2\omega+1}(2k) = 0, 
& f_{2\omega}(2k+1) = 0, & f_{2\omega+1}(2k+1) = f_{\omega}(k) \\
& g_{2\omega}(2k) = g_{\omega}(k), & g_{2\omega+1}(2k) = g_{\omega}(k),
& g_{2\omega}(2k+1) = h_{\omega}(k), & g_{2\omega+1}(2k+1) = f_{\omega}(k) \\
& h_{2\omega}(2k) = 0, & h_{2\omega+1}(2k) = g_{\omega}(k),
& h_{2\omega}(2k+1) = h_{\omega}(k), & h_{2\omega+1}(2k+1) = 0 \\
\end{array}
$$
\end{lemma}

\proof The proof is easy: it uses the definition of ${\omega \choose \ell}_2$, which 
in particular shows for any $2$-adic integer $\omega$ and any integer $\ell$ that
$$
\begin{array}{lll}
&\displaystyle{2\omega \choose 2\ell}_2 = 
\displaystyle{\omega \choose \ell}_2 {0 \choose 0}_2 = {\omega \choose \ell}_2  
\ \ \ \ \ \ \
&\displaystyle{2\omega+1 \choose 2\ell}_2 = 
\displaystyle{\omega \choose \ell}_2 {1 \choose 0}_2 = {\omega \choose \ell}_2 \\
&\displaystyle{2\omega \choose 2\ell+1}_2 = 
\displaystyle{\omega \choose \ell}_2 {0 \choose 1}_2 = 0 
\ \ \ \ \ \ \ 
&\displaystyle{2\omega+1 \choose 2\ell+1}_2 = 
\displaystyle{\omega \choose \ell}_2 {1 \choose 1}_2 = {\omega \choose \ell}_2.
\end{array}
$$

\begin{remark}
The sequences above occur in the OEIS \cite{OEIS} when $\omega = n$ is an integer. 
In particular, $({\frac{n+k}{2} \choose k})_{n,k} = ({n+k+1 \choose 2k+1})_{n,k}$ is equal 
to A168561; also $({n+k \choose 2k})_{n,k}$ is equal to A085478; finally, up to 
shifting $k$, we have that $({n+k \choose 2k+1})_{n,k}$ is equal to A078812.

We can also note that $f_{\omega}(k) \equiv g_{\omega}(k) + h_{\omega}(k) \bmod 2$,
for any integer $k \geq 0$.
\end{remark}

\begin{theorem}\label{algebr}

Suppose that $\lambda_n = 2^{n+1} - 1$. The following results then hold.

\medskip

-- The formal power series $Q_{\omega}(X)$ is either a polynomial 
if $\omega \in {\mathbb Z}$ or a transcendental series over ${\mathbb Q}(X)$ 
if $\omega \in {\mathbb Z}_2 \setminus {\mathbb Z}$. 

\medskip

-- The formal power series $Q_{\omega}(X)$ is algebraic over ${\mathbb Z}/2{\mathbb Z}(X)$
if and only if $\omega$ is rational. It is rational if and only if it is a polynomial, which
happens if and only if $\omega$ is a rational integer.

\end{theorem}

\proof The first assertion is a consequence of a classical theorem of Fatou \cite{Fatou} 
which states that {\it a power series $\sum_{n \geq 0} a_n z^n$ with integer 
coefficients that converges inside the unit disk is either rational or transcendental 
over ${\mathbb Q}(z)$}. This implies that the formal power series $Q_{\omega}(X)$ is 
either rational or transcendental over ${\mathbb Q}(X)$.
We then have to prove that if $Q_{\omega}$ is a rational
function, then it is a polynomial, or equivalently that $\omega$ is a rational integer 
(use Theorem~\ref{polynomial}). Now to say that $Q_{\omega}$ is rational is to say that the 
sequence of its coefficients is ultimately periodic, which implies that the sequence of their 
absolute values $(f_{\omega}(k))_{k \geq 0}=({\omega + k + 1 \choose 2k+1}_2)_{k \geq 0}$ 
is ultimately periodic. Let $\theta$ be its period. We have, for large $k$, that
${\omega + k + 1 \choose 2k+1}_2 = {\omega + k + \theta + 1 \choose 2(k+\theta)+1}_2$.
If $\theta$ is odd, the left side is zero for $\omega+k$ odd while the right side
is zero for $\omega+k$ even. Thus ${\omega + k + 1 \choose 2k+1}_2 = 0$ for large $k$,
and $Q_{\omega}$ is a polynomial. So suppose that $\theta$ is even. Let us suppose that
$\omega$ does not belong to ${\mathbb Z}$, then its $2$-adic expansion contains infinitely 
many blocks $01$. Consider the first such block: there exist $\ell \geq 0$ and $j \geq 0$ 
such that $\omega = 2^{\ell} - 1 + 2^{\ell +1}(2^j(2\omega'+1))$. Then for any integer $k'$ 
we have $f_{\omega}(2^{\ell} - 1 + 2^{\ell + 1}(2^j(2k'+1))) =
{\omega' + k' + 1 \choose 2k'+1}_2$. The sequence
$(f_{\omega}(2^{\ell} - 1 + 2^{\ell + 1}(2^j(2k'+1)))_{k' \geq 0}$ is ultimately
periodic and $\theta/2$ is a period. But from Lemma~\ref{digit} (vi) this sequence
is equal to $({\omega' + k' + 1 \choose 2k'+1}_2)_{k' \geq 0}$. As previously,
either $\theta/2$ is odd and this sequence is ultimately equal to zero or
$\theta/2$ is even. In the first case, as above, $\omega'$ belongs to ${\mathbb Z}$,
so does $\omega$, which is impossible. In the second case, we iterate the reasoning 
that used Lemma~\ref{digit} (vi), with $\omega$ replaced by $\omega'$ and $k$ by $k'$,
where the first block $01$ occurring in $\omega$ is replaced by the first such block
occurring in $\omega'$. The fact that $\theta$ cannot be divisible by arbitrarily
large powers of $2$ gives the desired contradiction.

\bigskip

In order to prove the second assertion, we first suppose that $Q_{\omega}(X)$ is 
algebraic over ${\mathbb Z}/2{\mathbb Z}(X)$. If $\omega=-1$, $Q_{\omega}(X) = 0$. 
Otherwise write $\omega = 2^{\ell} - 1 + 2^{\ell + 1}\omega'$ as in Lemma~\ref{digit} (v). 
The algebraicity of $Q_{\omega}(X)$ over ${\mathbb Z}/2{\mathbb Z}(X)$ implies that the 
sequence $({\frac{\omega+k}{2} \choose k}_2 \bmod 2)_{n \geq 0}$ is $2$-automatic 
(from a theorem of Christol, see \cite{Christol, CKMFR} or \cite{AS}). 
Using Lemma~\ref{kernel} (i) we thus have that the sequence
$({\omega+k+1 \choose 2k+1}_2)_{k \geq 0}$ is $2$-automatic. Thus its subsequence
obtained for $k = 2^{\ell} - 1 + 2^{\ell + 1} k'$, namely
$({\omega+ 2^{\ell} + 2^{\ell+1}k' \choose 2^{\ell+1}-1+2^{\ell+1}(2k')}_2)_{k' \geq 0}$
is also $2$-automatic (see, e.g., \cite[Theorem~6.8.1,~page~189]{AS}).
But this last sequence is equal to
$({2^{\ell+1} -1 + 2^{\ell+1}(\omega'+k') 
\choose 2^{\ell+1}-1+2^{\ell+1}(2k')}_2)_{k' \geq 0}$, i.e., to
$({\omega'+k' \choose 2k'}_2)_{k' \geq 0}$ (look at the $2$-adic expansions 
and use Definition~\ref{binom}). But this in turns implies (see, e.g., 
\cite[Corollary~5.5.3,~page~167]{AS}) that the subsequence
$({\omega'+2^j \choose 2^{j+1}}_2)_{j \geq 0}$ is ultimately periodic. Using
Lemma~\ref{digit} (ii) this means that $\omega$ is rational.

\medskip

Now suppose that $\omega$ is rational. Denote by $T\omega$ the $2$-adic integer
defined by $T\omega = (\omega - \omega_0)/2$ (i.e., $T\omega$ is the $2$-adic integer
obtained by shifting the sequence of digits of $\omega$). Also note $T^j$ the $j$-th
iteration of $T$. Define (with the notation of Lemma~\ref{kernel}) the set ${\cal K}$ by
$$
{\cal K} =: \left\{(f_{T^j\omega}(k))_{k \geq 0}, \ j \in {\mathbb N}\right\}
            \cup
            \left\{(g_{T^j\omega}(k))_{k \geq 0}, \ j \in {\mathbb N}\right\}
            \cup
            \left\{(h_{T^j\omega}(k))_{k \geq 0}, \ j \in {\mathbb N}\right\}.
$$
As a consequence of Lemma~\ref{kernel}, we see that ${\cal K}$ is stable by the 
maps defined on ${\cal K}$ by $(v_k)_{k \geq 0} \to (v_{2k})_{k \geq 0}$ and
$(v_k)_{k \geq 0} \to (v_{2k+1})_{k \geq 0}$ (use that for any $2$-adic integer
$\omega = \omega_0 \omega_1 \ldots$ one has $\omega = 2T\omega + \omega_0$). 
On the other hand we have from Lemma~\ref{digit}~(iv) that
${\frac{\omega+k}{2} \choose k}_2 = f_{\omega}(k)$.
Hence the $2$-kernel of the sequence $({\frac{\omega+k}{2} \choose k}_2)_{k \geq 0}$,
i.e., the smallest set of sequences containing that sequence and stable under the maps
$(v_k)_{k \geq 0} \to (v_{2k})_{k \geq 0}$ and
$(v_k)_{k \geq 0} \to (v_{2k+1})_{k \geq 0}$, is a subset of ${\cal K}$. Now, since $\omega$ 
is rational, the set of $2$-adic integers $\{T^j\omega, \ j \in {\mathbb N}\}$ is finite. Hence
the $2$-kernel of $({\frac{\omega+k}{2} \choose k}_2)_{k \geq 0}$ is finite and this sequence is 
$2$-automatic (see, e.g., \cite{AS}). This implies that the formal power series $Q_{\omega}(X)$
is algebraic over ${\mathbb Z}/2{\mathbb Z}(X)$ (using again Christol's theorem, see 
\cite{Christol, CKMFR} or \cite{AS}).

\medskip

Finally, $Q_{\omega}(X)$ reduced modulo $2$ is rational if and only if the sequence of its 
coefficients $(f_{\omega}(k))_{k \geq 0}=({\omega + k + 1 \choose 2k+1}_2)_{k \geq 0}$ 
modulo $2$ is ultimately periodic, which is the same as saying that the sequence
$(f_{\omega}(k))_{k \geq 0}=({\omega + k + 1 \choose 2k+1}_2)_{k \geq 0}$ itself is ultimately
periodic. But from the first part of the proof this implies that $Q_{\omega}(X)$ (not 
reduced modulo $2$) is a polynomial, hence that $Q_{\omega}(X)$ modulo $2$ is a polynomial. 
Conversely, if $Q_{\omega}(X)$ modulo $2$ is a polynomial, then the sequence of its coefficients 
$(f_{\omega}(k))_{k \geq 0}=({\omega + k + 1 \choose 2k+1}_2)_{k \geq 0}$ modulo $2$ is 
ultimately equal to $0$, and so is $(f_{\omega}(k))_{k \geq 0}$ not reduced modulo $2$. Thus
$Q_{\omega}(X)$ not reduced modulo $2$ is a polynomial, thus $\omega$ is a rational integer
by using Theorem~\ref{polynomial}.

\begin{remark}

\ { }

\begin{itemize}

\item The authors of \cite{AMFvdP} prove that the formal power series
$(1+X)^{\omega} = \sum_{k \geq 0} {\omega \choose k}_2 X^k$ is algebraic over 
${\mathbb Z}/2{\mathbb Z}(X)$ if and only if $\omega$ is rational. They do not ask when that 
series is rational, i.e., belongs to ${\mathbb Z}/2{\mathbb Z}(X)$, but this is clear since 
for rational $\omega = a/b$, with integers $a,b > 0$, we have 
$((1+X)^{\omega})^b \equiv (1+X)^a \bmod 2$. 
Hence if $(1+X)^{\omega}$ is a rational function A/B with $A$ and $B$ coprime polynomials, then
$A^b \equiv (1+X)^a B^b$ hence $B$ is constant, i.e., $(1+X)^{\omega}$ is a polynomial. Now if
$a < 0$ and $b > 0$, we have that $(1+X)^{-\omega}$ is a polynomial, hence $(1+X)^{\omega}$ is 
the inverse of a polynomial. Finally $(1+X)^{\omega}$ is a rational function if and only if
$\omega \in {\mathbb Z}$.

\item In the same vein, the authors of \cite{AMFvdP} prove that, if $\omega_1$, $\omega_2$,
..., $\omega_d$ are $2$-adic integers, then the formal power series $(1+X)^{\omega_1}$,
$(1+X)^{\omega_2}$, ..., $(1+X)^{\omega_d}$ are algebraically independent over 
${\mathbb Z}/2{\mathbb Z}(X)$ if and only if $1$, $\omega_1$, $\omega_2$,..., $\omega_d$ are
linearly independent over ${\mathbb Z}$. Is a similar statement true for the $Q_{\omega}$?

\item Another question is to ask whether a similar study can be done in the $p$-adic case 
(here $p = 2$). The two papers \cite{Carlitz65, Chan} might prove useful.

\item Results of transcendence, hypertranscendence, and algebraic independence of values for
the generating function of the Stern-Brocot sequence have been obtained very recently by
Bundschuh (see \cite{Bundschuh}, see also the references therein).

\item A last question is the arithmetic nature of the real numbers $A(\varepsilon, \omega, g)$ 
defined by $A(\varepsilon, \omega, g) = \sum_{k <\!\!< \frac{k+\omega}{2}} 
\sigma(k, \varepsilon) g^{-k}$ where $g \geq 2$ is an integer, the sequence $(\varepsilon_n)_n$ 
is ultimately periodic, and $\omega \in {\mathbb Z}_2 \setminus {\mathbb Z}$. Take in particular 
$\varepsilon = 0$ (thus $\sigma(k, \varepsilon) = (-1)^{\nu(k)}$). We already know that the number 
$A(0, \omega, g)$ is transcendental for $\omega \in ({\mathbb Q} \cap {\mathbb Z}_2) \setminus 
{\mathbb Z}$ by using \cite{AB}, the fact that the sequence $((-1)^{\nu(k)})_{k\geq 0}$ is 
$2$-automatic as recalled above, and the fact that the sequence 
$({\frac{k+\omega}{2} \choose k}_2)_{k \geq 0}$ is $2$-automatic for $\omega$ rational as 
seen in the course of the proof of Theorem~\ref{algebr} (the fact that $A(0, \omega, g)$ is 
not rational is a consequence of the non-ultimate periodicity of the sequence 
$((-1)^{\nu(k)}{\frac{k+\omega}{2} \choose k}_2)_{k \geq 0}$ for $\omega$ rational but not 
a rational integer, which has also been seen in the course of the proof of Theorem~\ref{algebr}).

\end{itemize}
\end{remark}

\end{document}